\documentclass[epic,eepic,epsfig,envcountresetchap,envcountsame,draft]{amsart} 
\usepackage{amssymb}
\usepackage{amsmath}
\usepackage{amscd}
\usepackage{epic,eepic,epsfig} 
\makeindex  

\newcommand{\Zb}{{\mathbb Z}}

\newcommand{\Gal}{{\rm Gal}}

\def\iii{{\Ic^\prime _2(X,Y)}}
 
\def\Cbb{{\mathcal C}}

\def\Nbb{{\mathcal N}}
\def\Pbb{{\mathcal P}}
\def\Qbb{{\mathcal Q}}

\def\Zbb{{\mathcal Z}}

\def\Cc{{\mathcal C}}

\def\Ic{{\mathcal I}}

\def\Mc{{\mathcal M}}

\def\Uc{{\mathcal U}}
\def\Vc{{\mathcal V}}

\def\01{{\overrightarrow{01}}}
\def\10{{\overrightarrow{10}}}

\def\hpb{\hfill $\Box$}

\def\lra{\longrightarrow}

\def\al{\alpha}
\def\be{\beta}
\def\de{\delta}
\def\De{\Delta}
\def\ga{\gamma}
\def\si{\sigma}
\def\La{\Lambda}
\def\Ga{\Gamma}

\def\mod{{\rm mod}}

\def\log{{\rm log}}

\def\pd1{{\partial \Delta [1]}} 
\def\d1{{\Delta [1]}}
\def\om{\omega}

\def\zl{\Zbb _\ell}

\def\zlt{\Zbb _\ell^\times}

\def\Cbb{{\mathbb C}}

\def\Nbb{{\mathbb N}}
\def\Pbb{{\mathbb P}}
\def\Qbb{{\mathbb Q}}

\def\Ybb{{\mathbb Y}}
\def\Zbb{{\mathbb Z}}

\def\Cc{{\mathcal C}}

\def\Ic{{\mathcal I}}

\def\Mc{{\mathcal M}}

\def\Uc{{\mathcal U}}
\def\Vc{{\mathcal V}}

\def\Zc{{\mathcal Z}}

\def\z2{{{\Zbb [{\frac{1}{2}}]}}}

\def\dfk{{\mathfrak d}}
\def\ffk{{\mathfrak f}}

\def\zfk{{\mathfrak z}}

\def\al{{\alpha}}
\def\be{{\beta}}
\def\ga{{\gamma}}
\def\si{{\sigma}}
\def\de{{\delta}}
\def\ka{{\kappa}}

\def\De{{\Delta}}

\def\zn{{z^{\frac{1}{\ell ^n}}}}

\def\01{{\overset{\to}{01}}}
\def\10{\overset{\to}{10}}

\def\hpb{\hfill $\Box$}

\def\lra{\longrightarrow}

\def\pk{ \pi _1(\Pbb ^1_{\bar K}\setminus \{ 0,1,\infty \} , \01 )}
\def\pq{  \Pbb ^1_{\bar \Qbb}\setminus \{ 0,1,\infty \}}
\def\qp{ \Qbb _p \{\{X,Y\}\}}
\def\mim{\mu ({\frac{i}{m}})}
\def\min{\mu ({\frac{i}{m}})}

\def\zp{\Zbb _p}
\def\zzp{\zp [[\zp]]}
\def\izd{\int _{\zp}x^{k-1}d}

\def\im{{\frac{i}{m}}}

\title[A polylogarithmic  measure associated with a path]{A polylogarithmic  measure associated with a path on $\Pbb ^1\setminus \{ 0,1,\infty \}$ and a $P$-adic Hurwitz zeta function}
\author{Zdzis{\l}aw Wojtkowiak}
\begin{document}

\date{\today}

\maketitle

\tableofcontents

\begin{abstract}
With every path on $\Pbb ^1_{\bar \Qbb}\setminus \{ 0,1,\infty \}$ there is associated a measure on $\Zbb _p$. 
The group $\Zbb _p^\times$ acts on measures. We consider two measures. One measure is associated to a path  from $\01$ to a root  of unity $\xi$ of order prime to $p$. 
Another measure is associated to a path  from $\01$ to  $\xi^{-1}$ and next it is acted by $-1\in \Zbb _p^\times$.
We show that the sum of these  measures  can be defined in a very elementary way.
Integrating against this sum of measures we get $p$-adic Hurwitz zeta functions constructed previously by Shiratani.
\end{abstract}

\section{Introduction}
\smallskip

Let $K$ be a number field, let $z\in \Pbb ^1(K)\setminus \{0,1,\infty \}$ and let $\ga$ be a path on $\Pbb ^1_{\bar K}\setminus \{ 0,1,\infty \}$ from $\01$ to $z$,
i.e. an isomorphism of the corresponding fiber functors. Let $p$ be a fixed prime number. The Galois group $G_K$ acts on
\[
 \pi _1(\Pbb ^1_{\bar K}\setminus \{ 0,1,\infty \},\01 )
\]
-- the pro-$p$ \'etale fundamental group. Let $\qp$ be the $\Qbb _p$-algebra of non-commutative formal power series in two non-commuting variables $X$ and $Y$. Let
\[
 E:\pk \to \qp
\]
be the continuous multiplicative embedding given by $E(x)=\exp X$ and $E(y)=\exp Y$, where $x$ and $y$ are standard generators of $\pk$. For any $\si \in G_K$ 
we define 
\[
 \ffk _\ga (\si ):=\ga ^{-1}\cdot \si (\ga )\in \pk
\]
and 
\[
 \Lambda _\ga (\si ):=E(\ffk _\ga (\si ))\in \qp \; .
\]
In the special case of the path $\pi$ from $\01$ to $\10$, the element $\ffk _\pi (\si )$ was studied by Ihara and his students (see \cite{I} and more other papers), 
Deligne 
(see \cite{D}), Grothendieck. The coefficients of the power series $\Lambda _\pi (\si)$ are analogues of the multi-zeta numbers studied already 
by Euler. 
For an arbitrary path $\ga$ the coefficients of the power series $\Lambda _\ga (\si )$ are analogues of values of iterated integrals evaluated at $z$. 

Observe that 
\[
 \Lambda _\ga (\si )\equiv 1+l_\ga (z)(\si )X\;\; {\rm modulo}\;\; I^2+(Y)
\]
for a certain $l_\ga (z)(\si )\in \Zb _p$, where $I$ is the augmentation ideal of $\qp$ and $(Y)$ is the principal ideal generated by $Y$. Let us set 
\[
 \Delta  _\ga (\si ):=\exp (- l_\ga (z)(\si )X)\cdot \Lambda _\ga (\si )\;.
\]
One possible way to calculate (some) coefficients of the power series $\Lambda _\pi(\si )$ and some other power series $\La _\ga (\si)$ is to use symmetries of
$\Pbb ^1 _{\bar \Qbb }\setminus \{0,1,\infty \}$, i.e. the so called Drinfeld-Ihara relations (see \cite{Dr} and \cite{I1}).  
For example in \cite{W7}, we have calculated even 
polylogarithmic coefficients of the power series $\La _\pi (\si )$ using the symmetries of $\Pbb ^1 _{\bar \Qbb }\setminus \{0,1,\infty \}$.

In \cite{NW} the authors have constructed a measure on $\Zbb _p$ for any path $\ga$ and expressed the $k$-th polylogarithmic coefficient of the power series 
$\log \De _\ga (\si )$ 
as integrals of the polynomial  $x^{k-1}$ against this measure recovering the old result of O. Gabber (see \cite{D0}). Let us denote this measure by $K(z)_\ga$.

\medskip

Now we shall describe the main result of this note. Let $m$ be a positive integer not divisible by $p$. Let us set 
\[
 \xi _m =\exp ({\frac{2 \pi \sqrt{-1}}{m}})\;.
\]
Let $0<i<m$. Further we chose paths $\be _i$ (resp. $\be _{m-i}$) on $\pq$ from $\01$ to $\xi _m^i$ (resp. $\xi _m^{m-i}$) such that $l_{\be _i}(\xi _m ^i)=0$ and 
$l_{\be _{m-i}}(\xi _m ^{m-i})=0$.

In \cite{W8} using the symmetry $\zfk  \mapsto 1/\zfk$ of $\pq$ we have shown that the polylogarithmic coefficient in degree $k$ of the formal power series  
\begin{equation}\label{eq:coefk}
 \log \Lambda _{\be _{m-i}}(\si )+(-1)^k \log \Lambda _{\be _{i}}(\si )
\end{equation}
is equal ${\frac{B_k({\frac{i}{m}})}{k!}}(1-\chi ^k (\si ))$, where $B_k(X)$ is the $k$-th Bernoulli polynomial and $\chi :G_{\Qbb (\mu _m)}\to \Zbb _p ^\times$
is the cyclotomic character (see \cite[Theorem 10.2.]{W8}). In this paper we shall calculate the same polylogarithmic coefficients using the measure
\[
 K(\xi _m^{m-i})_{\be _{m-i}}+\iota (K(\xi _m^i)_{\be _i})\, ,
\]
where $\iota$ is the complex conjugation acting on measures. To calculate these measures we use the symmetry  $\zfk  \to 1/\zfk$ of the tower of coverings
\[
 \Pbb ^1 _{\bar \Qbb}\setminus (\{0,\infty \}\cup \mu _{p^n})\to \pq ,\; \zfk \mapsto \zfk ^{p^n}
\]
of $\pq$.
However in contrast with the calculations in \cite{W8} we need to work only with terms in degree $1$. We show that the measure 
$K(\xi _m^{m-i})_{\be _{m-i}}+\iota (K(\xi _m^i)_{\be _i})$ is the sum of the Bernoulli measure $E_{1,\chi }$ 
(see \cite[the formula E.1 on page 38]{L})
and 
the measure we denote by $\mu _{\chi }({\frac{i}{m}})$. The definition of the measure $\mu _{\chi }({\frac{i}{m}})$ is 
very elementary and perhaps it is well known. From this it follows immediately the formula for the $k$-th polylogarithmic coefficient of the power series 
\eqref{eq:coefk}. The measure we got, allows to get the $p$-adic Hurwitz zeta functions as Mellin transform in the same way as the $p$-adic L-functions 
are the Mellin transforms of the measure $\psi E_{1,c}$, where $\psi $ is a character on $\Zbb _p^\times $ (see \cite[Chapter 4]{L}).

\bigskip


\section{An example of a measure on $\Zbb _p$}
\smallskip

This section can be seen as an attempt to construct a measure on $\Zbb _p$ which to a subset $a+p^n \Zbb _p$ associates $1/p^n$. We found the measure in question studying Galois actions on torsors of paths (see section 3). The measure is elementary and we think that it should be known.

\medskip
\noindent
If $a\in \Zbb _p$ and $a=\sum _{i=0}^\infty \al _i p^i$ with $0\leq \al _i\leq p-1$ then we set
\[
 v_n(a):=\sum _{i=0}^n \al _ip^i\;\;{\rm and}\;\;t_{n+1}(a):={\frac{a-v_n(a)}{p^{n+1}}}\, .
\]
Let us fix a positive integer $m>1$. For $k\in \Qbb ^\times $, $k={\frac{a}{b}}$ with $a,b\in \Zbb $ and $(b,m)=1$ we define 
\[
 [ k] _m\, \in \Nbb
\]
by the following two conditions
\[
 0\leq [ k] _m<m\;\;{\rm and}\;\; b  [ k] _m\equiv a \;\;{\rm modulo}\; m\,.
\]
Let us assume that $p$ does not divide $m$. Let $i$ be such that $0<i<m$. Observe that 
\begin{equation}\label{eq:p-np-m}
 [p^{-r}[ip^{-n}]_m   ]_m=[i  p^{-(n+r)}]_m\,.
\end{equation}
We define a sequence of integers 
\[
 (k_r(i))_{r\in \Nbb }
\]
by the equalities
\begin{equation}\label{eq:p-1p}
 p [ip^{-r}]_m   =[i  p^{-( r-1)}]_m +k_{r-1}(i)   m\, .
\end{equation}
Observe that 
\[
 0<{\frac{[i  p^{-( r-1)}]_m }{m}} <1 \;\;{\rm and}\;\;0<{\frac{p[ip^{-r}]_m }{m}}<p\,.
\]
Hence it follows that 
\[
 0\leq k_r(i)\leq p-1
\]
for all $r\geq 0$. Applying successively the formula \eqref{eq:p-1p} we get 
\begin{equation}\label{now}
 p^n [i  p^{-n}]_m  = i+\big(\sum _{\al =0}^{n-1} k_\al (i) p^\al \big)m\,. 
\end{equation}
It follows from \eqref{now} that 
\[
 {\frac{-i}{m}}=\sum _{\al =0}^{\infty } k_\al (i) p^\al
\]
and 
\begin{equation}\label{eq:develop}
  {\frac{i}{m}}=1+\sum _{\al =0}^{\infty }(p-1- k_\al (i)) p^\al\,.
\end{equation}
Another consequence of \eqref{now} is the equality
\[
 t_n(-{\frac{i}{m}})={\frac{-[ip^{-n}]_m}{m}}\,.
\]

For any positive integer $a$ such that $0\leq a <p^n$ we set 
\[
\de _n (a):=   \left\{
\begin{array}{ll}
-1 & \text{ if } a\geq 1+\sum _{\al =0}^{n-1 }(p-1- k_\al (i)) p^\al, \\
0 & \text{ if } a< 1+\sum _{\al =0}^{n-1 }(p-1- k_\al (i)) p^\al \;.
\end{array}
\right.
\]
\medskip
\noindent{\bf Definition-Proposition 1.1.} The function from the open-closed subsets of $\Zbb _p$ to $\Zbb _p$ defined by the formula
\[
 \mu ({\frac{i}{m}})(a+p^n\Zbb _p):={\frac{[ip^{-n}]_m}{m}}+\de _n (a)
\]
for $0\leq a<p^n$ is a measure.
\medskip

\noindent {\bf Proof.}  Let  $0\leq a<p^n$. We have
\[\sum _{b=0}^{p-1}\mu ({\frac{i}{m}})(a+bp^n+p^{n+1}\Zbb _p)= 
\sum _{b=0}^{p-1}({\frac{[ip^{-(n+1)}]_m}{m}}+\de _{n+1}(a+bp^n))=\]
\noindent
\[{\frac{p[ip^{-(n+1)}]_m}{m}}+\sum _{b=0}^{p-1}\de _{n+1}(a+bp^n)= 
{\frac{[ip^{-n }]_m}{m}}+k_n(i)+\sum _{b=0}^{p-1}\de _{n+1}(a+bp^n)\]
by the equality\eqref{eq:p-1p}. Observe that 
\[
\sum _{b=0}^{p-1}\de _{n+1} (a+bp^n):=   \left\{
\begin{array}{ll}
-k_n(i)-1 & \text{ if } a\geq 1+\sum _{\al =0}^{n-1 }(p-1- k_\al (i)) p^\al, \\
-k_n(i) & \text{ if } a< 1+\sum _{\al =0}^{n-1 }(p-1- k_\al (i)) p^\al \;.
\end{array}
\right.
\]
Hence finally we get 
$ \sum _{b=0}^{p-1}\mim (a +bp +p^{n+1}\Zbb _p)={\frac{[ip^{-n }]_m}{m}}+\de _n(a)=\mim (a+p^n\Zbb _p)$. \hpb

\medskip
\noindent{\bf  Proposition 1.2.} For $k\geq 1$ we have
\begin{enumerate}
 \item [i)] 
\[
 \int _{\zp}x^{k-1}d\mim (x)={\frac{1}{k}}\big( B_k({\frac{i}{m}})-B_k\big)\,,
\]
 \item[ii)] 
 \[
 \int _{\zp ^\times } x^{k-1} d \min (x)={\frac{1}{k}}\big(B_k( {\frac{i}{m}} )-B_k\big)-
 {\frac{p^ {k-1}}{k}}\big(B_k( {\frac{[ip^ {-1}]_m}{m}} )-B_k\big)\,.
 \]
\end{enumerate}

\medskip
\noindent {\bf Proof.} First we shall prove the formula i). Let us calculate the Riemann sum
\[
 \sum _{\al =0}^{p^n-1} \al ^{k-1} \min (\al +p^n\zp ) = \sum _{\al =0}^{p^n-1} \al ^{k-1}\big( {\frac{[ip^{-n}]_m}{m}}+\de _n (\al )\big)=
\]
\[
{\frac{[ip^{-n}]_m}{m}} \sum _{\al =0}^{p^n-1} \al ^{k-1}-\sum _{\al =0}^{p^n-1} \al ^{k-1}+\sum _{\al =0}^{v_{n-1}({\frac{i}{m}})-1} \al ^{k-1}\,.
\]
Observe that 
$$\sum _{\al =0}^{v_{n-1}({\frac{i}{m}})-1 } \al ^{k-1}={\frac{1}{k}}\big(B_k(v_{n-1}({\frac{i}{m}}))-B_k\big)$$
and it tends to ${\frac{1}{k}}\big(B_k( {\frac{i}{m}} )-B_k\big)$ if $n$ tends to $\infty$.  Hence the formula i) of the proposition follows because 
$ \sum _{\al =0}^{p^n-1} \al ^{k-1}$ tends to $0$ if $n$ tends to $\infty$\,.   

Observe that 
$$ \int _{\zp ^\times } x^{k-1} d \min (x)= \izd \min (x)-\int _{p\zp   }  x^{k-1}   d \min (x)\,.$$
We shall calculate Riemann sums for the integral $\int _{p\zp   }  x^{k-1}   d \min (x)$. We have
\[
\sum _{\al =0}^{p^n-1}(p\al)^{k-1}\mu (p\al +p^{n+1}\zp)=\sum _{\al =0}^{p^n-1}p ^{k-1}\al ^{k-1}{\frac{[ip^{-(n+1)}]_m}{m}} 
 +\sum _{\al =0}^{p^n-1}p ^{k-1}\al ^{k-1}\delta _{n+1}(p\al)\, .
\]
The first sum tend to $0$ if $n$ tends to $\infty$. Observe that 
\[
 \sum _{\al =0}^{p^n-1}p ^{k-1}\al ^{k-1}\delta _{n+1}(p\al)=\sum_{0<\al<p^n,\, p\al\geq v_n({\frac{i}{m}})}p ^{k-1}\al ^{k-1}(-1)=
\]
\[
 -\sum _{\al =0}^{p^n-1}p ^{k-1}\al ^{k-1}+\sum_{0<\al<p^n,\, p\al< v_n({\frac{i}{m}})}p ^{k-1}\al ^{k-1}\,.
\]
Let $0\leq \be _0 <p$ be such that $v_n({\frac{i}{m}})\equiv \be _0$ modulo $p$.
Then
\[
v_{n-1}({\frac{[ip^{-1}]_m}{m}}) =   \left\{
\begin{array}{ll}
1+{\frac{1}{p}}(v_n({\frac{i}{m}})-\be _0) & \text{ if } \be _0\neq 0, \\
{\frac{1}{p}}v_n({\frac{i}{m}}) & \text{ if } \be _0=0 \;.
\end{array}
\right.
\]
Hence it follows that
\[\sum_{0<\al<p^n,\, p\al< v_n({\frac{i}{m}})}p ^{k-1}\al ^{k-1}=p^{k-1}\sum _{\al=0}^{v_{n-1}({\frac{[ip^{-1}]_m}{m}})-1}\al ^{k-1}.\]
If $n$ tends to $\infty$ the last sum tends to $p^{k-1}{\frac{1}{k}}\big(B_k ({\frac{[ip^{-1}]_m}{m}})-B_k\big)\,.$
Hence the proof of the formula ii) is finished. \hpb

\medskip

If $c\in \zp ^\times \setminus \mu _{p-1}$ we define 
\begin{equation} \label{eq:1.6.}
 \mu _c ({\frac{i}{m}}):=\min -c\min \circ c^{-1}\,. 
\end{equation}
Then we have
\begin{equation} \label{eq:1.7.}
 {\frac{1}{1-c^k}}\izd \mu _c ({\frac{i}{m}})(x)= {\frac{1}{k}}\big(B_k( {\frac{i}{m}} )-B_k\big)\,. 
\end{equation}

\medskip

\noindent{\bf Corollary 1.3.} 
 Let  $P:\zp [[\zp ]]\to \zp [[T]]$ be the Iwasawa isomorphism given by $P(1)=1+T$. Then
 \[
  P(\min )(T)={ \frac{(1+T)^{{\frac{i}{m}}}-1}{T}}
 \]
and
\[
 P(\mu _c ({\frac{i}{m}}))={\frac{(1+T)^{{\frac{i}{m}}}-1}{T}} - {\frac{c\big( (1+T)^{c{\frac{i}{m}}}-1\big)}{(1+T)^c-1}}\,.
\]

\medskip
\noindent {\bf Proof.} The power series $ P(\min )(\exp X-1)$ is equal $\sum _{k=0}^\infty \big( \int _{\Zbb _p}x^k d\min (x)\big) X^k$.
Hence by the point i) of Proposition 1.2 it is equal 
\[
\sum _{k=0}^\infty {\frac{1}{(k+1)!}}\big(B_{k+1} ({\frac{i}{m}})-B_{k+1}\big)X^k\,.
\]
It follows from the definition of the Bernoulli numbers and the Bernoulli polynomials that this power series is equal ${\frac{\exp {\frac{i}{m}}X-1}{\exp X-1}}$.
Replacing $\exp X$ by $1+T$ we get the power series $P(\min )(T)$. \hpb

\medskip

We denote by 
\[
 \om :\zp ^\times \to \mu _{p-1}\subset \zp ^\times
\]
the Teichm\"uller character. For $x\in \zp ^\times $ we set 
\[
 [x]:=x\om (x)^{-1}\,.
\]
Let us define 
\[
 {\tilde H}_p(1-s, \om ^b,\im ):=\int _{\zp ^\times }[x]^s x^{-1}\om (x)^bd \min (x)\,.
\]

\medskip

\noindent{\bf Proposition 1.4.}  Let $k\equiv b$ modulo $p-1$. Then
\[
 {\tilde H}_p(1-k, \om ^b,\im )={\frac{1}{k}}\big(B_k( {\frac{i}{m}} )-B_k\big)-
 {\frac{p^ {k-1}}{k}}\big(B_k( {\frac{[ip^ {-1}]_m}{m}} )-B_k\big)\,.
\]

\medskip

\noindent {\bf Proof.} We have 
\[
  {\tilde H}_p(1-k, \om ^b,\im )=\int _{\zp ^\times }[x]^k x^{-1}\om (x)^bd \min (x)= \int _{\zp ^\times }x^{k-1} d \min (x)\,.
\]
Hence the proposition follows from the formula ii) of Proposition 1.2. \hpb

\medskip

\noindent
{\bf Remark 1.5.} A function closely related to our function $\tilde H_p(1-s,\om ^b,{\frac{i}{m}})$ appears in a paper of Shiratani 
(see \cite[Theorem 1, case $p\nmid f$]{Sh}).
 
\bigskip


\section{Action of the complex conjugation on measures}
\smallskip
We define an action of $\zp ^\times $ on the group ring $\zp [\zp ]$ by the formula
\[
 \al (\sum _{i=1}^na_i(x_i))=\al \sum _{i=1}^na_i(\al^{-1}x_i)
\]
and we extend by continuity to the action of $\zp ^\times$ on $\zp [[\zp]]$. The action of $\,-1\in \zp ^\times$ we denote by $\iota$. 
Then
\[
 \zzp=\zzp ^+\oplus \zzp ^-\, ,
\]
where $\iota$ acts on $\zzp ^+$ (resp. on $\zzp ^-$) as the identity (resp. as the multiplication by $-1$). For any $\mu \in \zzp$ we have the decomposition
\[
 \mu =\mu ^++\mu ^-\, ,
\]
where $\mu ^+={\frac{1}{2}}(\mu +\iota (\mu ))\in \zzp ^+$ and $\mu ^-={\frac{1}{2}}(\mu -\iota (\mu ))\in \zzp ^-$. Observe that
\begin{equation}\label{eq:2.1.}
 \izd \iota (\mu)=(-1)^k\izd \mu\, .
\end{equation}
Hence it follows 
\begin{equation}\label{eq:2.2.}
 \izd \mu ^+:=   \left\{
\begin{array}{ll}
0 & \text{ for } \, k\; \text{odd}, \\
\izd \mu  & \text{ for } \, k\; \text{even}
\end{array}
\right.
\end{equation}
and
\begin{equation}\label{eq:2.3.}
 \izd \mu ^-:=   \left\{
\begin{array}{ll}
 \izd \mu & \text{ for } \, k\; \text{odd}, \\
0  & \text{ for } \, k\; \text{even}\,.
\end{array}
\right.
\end{equation}

\medskip

In \cite[Proposition 10.5]{W8} we have shown that
\begin{equation}\label{eq:2.4.}
 \izd (K(\xi _m^{-i})+K(\xi _m^i))={\frac{1}{k}}B_k({\frac{i}{m}})(1-\chi ^k)\;\text{for}\;k\;\text{even}
\end{equation}
and
\begin{equation}\label{eq:2.5.}
 \izd (K(\xi _m^{-i})-K(\xi _m^i))={\frac{1}{k}}B_k({\frac{i}{m}})(1-\chi ^k)\;\text{for}\;k\;\text{odd}.
\end{equation}
\medskip
Hence it follows from \eqref{eq:2.2.} and \eqref{eq:2.3.} that  
\begin{equation}\label{eq:2.6.}
 \izd \big( (K(\xi _m^{-i})+K(\xi _m^i))^+ + (K(\xi _m^{-i})-K(\xi _m^i))^-\big)= 
\end{equation}
\[
 {\frac{1}{k}}B_k({\frac{i}{m}})(1-\chi ^k)\;\text{for}\;k\geq 1\,.
\]

Observe that 
\[
 (K(\xi _m^{-i})+K(\xi _m^i))^+ + (K(\xi _m^{-i})-K(\xi _m^i))^-=K(\xi _m^{-i})+\iota (K(\xi _m^i)) \,.
\]
Hence we get
\begin{equation}\label{eq:2.7.}
 \izd (K(\xi _m^{-i})+\iota (K(\xi _m^i)))={\frac{1}{k}}B_k({\frac{i}{m}})(1-\chi ^k)\;\;\text{for}\;\;k\geq 1\,.
\end{equation}

\medskip
The proof of the formulas \eqref{eq:2.4.} and \eqref{eq:2.5.} given in \cite{W8} is based on the symmetry $\zfk \mapsto 1/\zfk$ of $\pq$ and the study of the polylogarithmic coefficients 
( at $YX^{k-1}$) of the power series $\Lambda _{\be _i}(\si ) $ and  $\Lambda _{\be _{m-i}}(\si ) $.
Recently, H. Nakamura (see \cite{N}) got these formulas using directly the inversion formula from \cite[section 6.3]{NW2}.  

In this paper we calculate explicitely the measure $K(\xi _m^{-i})+\iota (K(\xi _m^i))$. We use also the symmetry $\zfk \mapsto 1/\zfk$ of the tower of coverings
\[
 \Pbb 1_{\bar \Qbb}\setminus (\{0,\infty\}\cup \mu _{p^n})\to \pq,\; \zfk \mapsto \zfk ^{p^n}
\]
but only in degree $1$.

The third possible method to calculate the measure $K(\xi _m^{-i})+\iota (K(\xi _m^i))$ is to use the explicit formula for measures $K(z)$ 
(see \cite[Proposition 3]{NW}). Compare the three different proofs of Proposition 5.13 in \cite{NW2}. Two proofs are given in \cite{NW2} and the third one in \cite{W8}
(the second proof of Lemma 4.1.)

\bigskip


\section{Measures associated with roots of unity}
\smallskip
We set 
\[
 \xi _r:=\exp ({\frac{2\pi \sqrt{-1}}{r}})
\]
for a natural number $r$. 
Let us set
\[
 V_n:=\Pbb ^1_{\bar \Qbb}\setminus (\{0,\infty\}\cup \mu _{p^n}).
\]
We recall that $\pi _1(V_n,\01 )$ - pro-$p$  \'etale fundamental group - is free on generators $x_n$ - loop around $0$ - 
and $y_{n,i}$ - loops around $\xi _{p^n}^i$   for $0\leq i<p^n$.
 
\medskip
 
For each $0<i<m$, let $\al _i$ be a path on $V_0=\pq$ from $\01$ to $\xi _m^i$ which is the composition of an arc from $\01$ to ${\overset{\longrightarrow}{0\xi _m^i}}$
in an infinitesimal neighbourhood of $0$ followed by the canonical path (straight line) from ${\overset{\longrightarrow}{0\xi _m^i}}$ to $\xi _m^i$.

\medskip

Let us set
\[
 \be _i:=\al _i \cdot x^ {-{\frac{i}{m}}}\,.
\]
Observe that $l(\xi _m^i)_{\be _i}=0$. If we regard the path $\al _i$ as the path on $V_n$ then we denote it by 
\[
 \, _n\al _i\,.
\]
Then
\[
 \, _n \be _i:= \, _n\al _i \cdot x_n^ {-{\frac{i}{m}}}
\]
is also a path on $V_n$. Let
\[
 \tilde \be _i ^n \;\;(\text{resp.}\;\; \tilde \al _i^n \, )
\]
be the lifting of $\be _i$ (resp. $\al _i$) to $V_n$ starting from $\01$.
\medskip
Let $0\leq j<p^n$. We denote by $s_n^j$ a lifting of $x_0^j$ to $V_n$ starting from $\01$. Observe that $s_n^j$ is a path on $V_n$ from $\01$ to 
${\overset{\longrightarrow}{0\xi _{p^n}^j}}$

\medskip
\noindent{\bf Lemma 3.1.} We have
\[
 \tilde \be _i^n=\, _n\be _{[ip^{-n}]_m}=\, _n\al _{[ip^{-n}]_m}\cdot x_n^{-{\frac{[ip^{-n}]_m}{m}}}\,.
\]

\medskip

\noindent {\bf Proof.} Observe that the lifting of $x^ {-{\frac{i}{m}}} $ to $V_n$ is equal $ s_n ^{v_{n-1}(-{\frac{i}{m}})}\cdot x_n^{t_n(-{\frac{i}{m}})}$.
The lifting of $\al _i$ to $V_n$ is a path (an arc) from $\01$ to ${\overset{\longrightarrow}{w}} 
:={ \overset{\longrightarrow} {0 \xi _{p^nm}^i }}$ in the positive sense composed 
with the canonical path from ${\overset{\longrightarrow}{w}}$ to $\xi_{p^n m}^i$. 
Hence the  lifting of $\be _i$ is the composition of $ s_n ^{v_{n-1}(-{\frac{i}{m}})}\cdot x_n^{t_n(-{\frac{i}{m}})}$ with the lifting of $\al _i$ multiplied by 
$\xi _{p^n}^{v_{n-1}(-{\frac{i}{m}})}$.

We have 
\[
\xi _{p^n}^{ v_{n-1}(-{\frac{i}{m}}) }  \xi _{p^n m}^i=\xi _{p^n m}^{ mv_{n-1}(-{\frac{i}{m}}) +i} \,.
\] 
 Observe that $0\leq v_{n-1}(-{\frac{i}{m}})\cdot m+i<p^n m$ and that $p^n$ divides $v_{n-1}(-{\frac{i}{m}})\cdot m+i$. Moreover we have
 ${\frac{v_{n-1}(-{\frac{i}{m}})\cdot m+i}{p^n}}\cdot p^n\equiv i$ modulo $m$. Hence it follows that ${\frac{v_{n-1}(-{\frac{i}{m}})\cdot m+i}{p^n}}=[ip^{-n}]_m$. 
 Therefore we get 
 \[
 -{\frac{[ip^{-n}]_m}{m}}=-{\frac{1}{p ^n}}(v_{n-1}(-{\frac{i}{m}}) +{\frac{i}{m}}) =t_n( -{\frac{i}{m}} ).  
 \]
 Hence it
follows that the lifting of $\be _i$ is $\, _n\al _{[ip^{-n}]_m}\cdot x_n^{-{\frac{[ip^{-n}]_m}{m}}}$ . \hpb

\medskip

To simplify the notation we set
\[
 r_n =[ip^{-n}]_m\;\;\text{and}\;\;v_{n-1}= v_{n-1}(-{\frac{i}{m}})\,.
\]
Then we have
\[
  \tilde \be _i^n=\, _n\al _{r_n}\cdot x_n^{-{\frac{r_n}{m}}}\;\;\text{and}\;\; \tilde \be _{m-i}^n=\, _n\al _{m-r_n}\cdot x_n^{{\frac{r_n}{m}}-1}
\]

Let $h:V_n \to V_n$ be given $\zfk \to 1/\zfk$. Let $p_n$ be the canonical path from $\01$ to $\10$ on $V_n$, $t_n$ a path from $\10$  to ${\overset{\to}{ 1\infty}}$
(half circle in the positive sense in an infinitesimal neighbourhood of $1$) and $q_n=h(p_n)$. We set
\[
 \Gamma _n:=q_n\cdot t_n\cdot p_n\,.
\]

\bigskip

\noindent{\bf Lemma 3.2.} We have
\[
 \tilde \be _{m-i}^n=h(\tilde \be _{i}^n )\cdot \Ga _n \cdot z_n^{\frac{r_n}{m}}\cdot x_n\cdot y_{n,-1}\cdot\ldots y_{n,-v_{n-1}}\cdot x_n ^{{\frac{r_n}{m}}-1}
\]
in $\pi _1(V_n,\01 )$.

\medskip

\noindent{\bf Proof.}  One checks that $\, _n\al _{m-r_n}=h(\, _n\al _{ r_n})\cdot \Ga _n \cdot x_n\cdot y_{n,-1}\cdot\ldots y_{n,-v_{n-1}}$.
The formula of the lemma follows from Lemma 3.1.
\hpb

\bigskip

\noindent{\bf Lemma 3.3.} Let $\si \in G_{\Qbb (\mu _m)}$. Then writting additively we have
\[
 \ffk _{\Gamma _n}(\si )\equiv \sum _{k=0}^{p^n-1}E^{(n)}_{1,\chi (\si )}(k)y_{n,k}\;\;\text{modulo}\;\;(\pi _1(V_n,\01),\pi _1(V_n,\01 ))\;.
\]

\medskip

\noindent{\bf Proof.} See the proof of Lemma 4.1 in \cite{W8} or the second proof of Proposition 5.13 in \cite{NW2}. \hpb

\medskip

It follows from Lemma 3.2 that 
\[
  \ffk _{\tilde \be _{m-i}^n}(\si ) \equiv \Gamma _n^{-1}h(\ffk _{\tilde \be _{ i}^n}(\si ))\Ga _n \cdot \ffk _{\Ga _n}(\si) \cdot   
\]
\[
 \big(z_n ^{\frac{r_n}{m}}\cdot x_n\cdot y_{n,-1}\cdot\ldots y_{n,-v_{n-1}}\cdot x_n ^{{\frac{r_n}{m}}-1}\big) ^{-1}\cdot 
 \si \big(z_n ^{\frac{r_n}{m}}\cdot x_n\cdot y_{n,-1}\cdot\ldots y_{n,-v_{n-1}}\cdot x_n ^{{\frac{r_n}{m}}-1}\big)
\]
modulo $(\pi _1(V_n,\01),\pi _1(V_n,\01 ))$. Hence writting the result additively we get
\[
 \sum _{k=0}^{p^n-1}K^{(n)}(\xi _m ^{-i})(\si )(k)y_{n,k}\equiv \sum _{k=0}^{p^n-1}K^{(n)}(\xi _m ^{ i})(\si )(k)y_{n,-k}+
 \sum _{k=0}^{p^n-1}E^{(n)}_{1,\chi (\si )}(k)y_{n,k}+
\]
\[
 \sum _{k=0}^{p^n-1}(1-\chi (\si )) { \frac{[ip^{-n}]_m}{m} } y_{n,k}-\sum _{j=1}^ { v_{n-1}(- {\frac{i}{m}} )}  y_{n,-j}+
 \chi (\si )\sum _ {j=1}^{ v_{n-1}(-{ \frac{i}{m} })}y_{n,-[j\chi (\si )]_ {p^n}} 
\]
modulo $(\pi _1(V_n,\01),\pi _1(V_n,\01 ))$. Observe that  $ v_{n-1}( {\frac{i}{m}}) =p ^n- v_{n-1}(-{\frac{i}{m}}) $, 
Hence the last two sums we can rewrite in the form
\[
 \sum _{j= v_{n-1}( { \frac{i}{m} })}^ {p^n-1}y_{n,j}+\chi (\si )\sum _{j= v_{n-1}( {\frac{i}{m}})}^{p^n-1}y_{n,[j\chi (\si )]_{p^n}}\,.
\]

Comparing coefficients at $y_{n,k}$ we get for $0\leq k<p^n$ 
\begin{equation}\label{eq:3.4.}
K^{(n)}(\xi _m^{-i})(\si )(k)- K^{(n)}(\xi _m^{i})(\si )(-k)=
\end{equation}
\[
E^{(n)}_{1,\chi (\si)}(k)+{\frac{[ip^{-n}]_m}{m}}+\de _n(k)-\chi (\si ){\frac{[ip^{-n}]_m}{m}}
+\chi (\si )\de _n([\chi (\si )^{-1}k]_{p^n})=
\]
\[
E^{(n)}_{1,\chi (\si)}(k)+\mu _{\chi (\si )}({\frac{i}{m}})(k)
\]

by the definition of the measure $\mu _{\chi (\si )}({\frac{i}{m}})$.

\medskip

\noindent
{\bf Theorem 3.5.}   Let $m$ be a positive integer not divisible by $p$ and let $0<i<m$. Then we have
\[
 K(\xi _m ^{-i})(\si )+\iota (K(\xi _m^i)(\si ))=E_{1,\chi (\si )}+\mu _{\chi (\si )}({\frac{i}{m}})\,.
\]

\medskip

\noindent {\bf Proof.} The theorem follows from the formula \eqref{eq:3.4.}. \hpb

\medskip

\noindent
{\bf Corollary 3.6.}   Let $\si \in G_{\Qbb (\mu _m)}$ be such that $\chi (\si )^{p-1}\neq 1$.   Then we have
\begin{enumerate}
 \item [i)] 
\[
 {\frac{1}{1-\chi (\si )^k}} \izd \big( K(\xi _m ^{-i})(\si )+\iota (K(\xi _m^i)(\si ))\big)={\frac{B_k({\frac{i}{m}})}{k}}\,,
\]
 \item [ii)] 
 \[
  P(K(\xi _m ^{-i})(\si )+\iota (K(\xi _m^i)(\si )))(T)={ \frac{ (1+T)^{\frac{i}{m}} }{T} }-  
  { \frac { \chi (\si)(1+T)^ {\chi (\si ){\frac{i}{m}}}} {(1+T)^{\chi (\si )}-1} }\,.
 \]
\end{enumerate}

\medskip

\noindent {\bf Proof.} The point i) of the corollary follows from Theorem 3.5 and the formula \eqref{eq:1.7.}.  The point ii) follows immediately from Corollary 1.8 
and the equality 
$P(E_{1,\chi (\si )})(T)={\frac{1}{T}}-{\frac{\chi (\si)}{(1+T)^{\chi (\si )}-1}}$. \hpb

\medskip
Now we define 
\[
 L^\be (1-s,(\xi _m^{-i})+\iota (\xi _m ^i);\si ):=
\]
\[
 {\frac{1}{1-\om (\chi (\si ))^\be [\chi (\si )]^s }}\int _{\Zbb ^\times _p}[x]^s x^{-1}\om (x)^\be d \big( (K(\xi _m^{-i})+\iota ( K(\xi _m^i))\big) (\si )\,.
\]

\medskip

\noindent
{\bf Theorem 3.7.} Let $\si \in G_{\Qbb(\mu _m)}$ be such that $\chi (\si )^{p-1}\neq 1$.
\begin{enumerate}
 \item [i)] Let $k\equiv \be $ modulo $ (p-1)$. Then
\[
   L^\be (1-k,(\xi _m^{-i})+\iota (\xi _m ^i);\si )={\frac{1}{k}}B_k({\frac{i}{m}})-p^{k-1}{\frac{1}{k}}B_k({\frac{[ip^{-1}]_m}{m}})\,.
\]
 \item [ii)] Let $\si , \si _1 \in G_{\Qbb(\mu _m)}$ be such that $\chi (\si )^{p-1}\neq 1$ and $\chi (\si _1 )^{p-1}\neq 1$.Then
 \[
 L^\be (1-s,(\xi _m^{-i})+\iota (\xi _m ^i);\si )=L^\be (1-s,(\xi _m^{-i})+\iota (\xi _m ^i);\si _1)\,,
 \]
 i.e. the function $L^\be (1-s,(\xi _m^{-i})+\iota (\xi _m ^i);\si )$ does not depend on $\si $.
\end{enumerate}

\medskip

\noindent{\bf Proof.}  For $ k\equiv \be $ modulo $p-1$ we have
\[
 L^\be (1-k,(\xi _m^{-i})+\iota (\xi _m ^i);\si )={\frac{1}{1-\chi (\si )^k }}  \int _{\Zbb ^\times _p} x^{k-1}d(\mu _{\chi (\si )}({\frac{i}{m}})+E_{1,\chi (\si )})
\]
by Theorem 3.5. It follows from \cite[Theorem 2.3]{L} that 
${\frac{1}{\chi (\si )^k-1}}\izd E_{1,\chi (\si)}=-{\frac{1}{k}}B_k$. The ``periodicity"property $E_{1,\chi (\si)}^{(n)}(i)=E_{1,\chi (\si)}^{(n+1)}(pi)$ 
of the measure $E_{1,\chi (\si)}$ implies that 
\begin{equation}\label{EQ0}
  {\frac{1}{1-\chi (\si )^k}} \int _{\Zbb ^\times _p} x^{k-1}d E_{1,\chi (\si )}=(1-p^{k-1}){\frac{1}{k}}B_k\,.
\end{equation}

\medskip

Integrating the function $x^{k-1}$ against the measure $\mu _{\chi (\si )}({\frac{i}{m}})$ we get
\[
{\frac{1}{\chi (\si )^k-1}} \int _{\Zbb ^\times _p} x^{k-1}d \mu _{\chi (\si )}({\frac{i}{m}})(x)=
\]
\[
{\frac{1}{\chi (\si )^k-1}} \Big( \int _{\Zbb ^\times _p} x^{k-1}d\mu  ({\frac{i}{m}})(x) -
\int _{\Zbb ^\times _p} x^{k-1}d(\chi (\si )\mu  ({\frac{i}{m}})\circ \chi (\si )^{-1})(x)\Big)\,.
\] 
Observe that 
$\int _{\Zbb ^\times _p} x^{k-1}d(\chi (\si )\mu  ({\frac{i}{m}})\circ \chi (\si )^{-1})(x)=\chi (\si ) ^k
\int _{\Zbb ^\times _p} y^{k-1}d  \mu  ({\frac{i}{m}})(y)$
if we set $\chi (\si )y=x$. It follows from Proposition 1.9 that 
\begin{equation}\label{EQ1}
  {\frac{1}{\chi (\si )^k-1} }  \int _{\Zbb ^\times _p} x^{k-1}d\mu _  { \chi (\si )}( {\frac{i}{m}} )=
  {\frac{1}{k}}  
  \big(B_k( {\frac{i}{m}} )
  -B_k\big) -p^{k-1} 
  {\frac{1}{k}} 
  \big(B_k(  {\frac{[ip^{-1}]_m}{m}}   )-B_k\big)\,.
\end{equation}
After the addition of \eqref{EQ0} and \eqref{EQ1} we get the point i) of the theorem. 

Concerning the point ii) observe that the functions 
$L^\be (1-s,(\xi _m^{-i})+\iota (\xi _m ^i);\si )$ and $L^\be (1-s,(\xi _m^{-i})+\iota (\xi _m ^i);\si _1 )$ coincide for $k\equiv \be $ modulo $(p-1)$. 
Hence these functions are equal because they are equal on a dense subset of $\Zbb _p$. \hpb

\bigskip

\vglue 2cm

\vglue 1cm

\noindent Universit\'e de Nice-Sophia Antipolis

\noindent D\'epartement de Math\'ematiques

\noindent Laboratoire Jean Alexandre Dieudonn\'e

\noindent U.R.A. au C.N.R.S., N$^{\rm  o}$ 168

\noindent Parc Valrose -- B.P. N$^{\rm  o}$ 71

\noindent 06108 Nice Cedex 2, France

\smallskip

\noindent {\it E-mail address} wojtkow@math.unice.fr

\noindent {\it Fax number} 04 93 51 79 74

\medskip


\begin{thebibliography}{999}
\bibitem{D} {\sc P. Deligne}, Le groupe fondamental de la droite projective moins trois points, {\it  in} Galois Groups over Q (ed. Y.Ihara, K.Ribet and J.-P. Serre),
{\it Mathematical Sciences Research Institute Publications}, {\bf 16} (1989), pp. 79-297.
\bibitem{D0} {\sc P. Deligne}, letter to Grothendieck, 19.11.82.
\bibitem{Dr} {\sc V. Drinfeld}, On quasi-triangulated quasi-Hopf algebras and some groups closely associated with $\Gal (\bar \Qbb /\Qbb )$, 
Algebra: Analiz, 2(1990), pp. 114-148.
\bibitem{I}{\sc Y. Ihara},  {Profinite braid groups, Galois representations and complex multiplications},
Annals of Math. 123 (1986), pp. 43-106.
\bibitem{I1} {\sc Y. Ihara}, Braids, Galois Groups and Some Arithmetic Functions, Proc. of the Int. Congress of Math.  Kyoto 1990, Springer-Verlag pp. 99-120.
\bibitem{L} {\sc S. Lang}, Cyclotomic fields I and II, Springer-Verlag New York Inc. 1990.
\bibitem{N} {\sc H. Nakamura}, e-mail letter, October 27, 2014.
\bibitem{NW} {\sc H. Nakamura, Z. Wojtkowiak}, On the explicit formulae for $l$-adic polylogarithms, {\it in} Arithmetic Fundamental Groups and Noncommutative Algebra, {\it Proc. of Symposia in Pure Math.} {\bf 70}, AMS 2002, pp. 285-294.
\bibitem{NW2} {\sc H. Nakamura, Z. Wojtkowiak},  Homotopy and tensor conditions for  functional equations of $l$-adic and classical iterated integrals,
{\it in} in Non-abelian Fundamental Groups and Iwasawa Theory, London Math. Soc, Lecture Note Series, 393
pages 258--310,  2012, Cambridge UP.  
\bibitem{Sh} {\sc K. Shiratani}, On a Kind of p-adic Zeta Functions, {\it  in} Algebraic Number Theory (ed. S. Iyanaga), International Symposium, Kyoto 1976, pp. 213-217.

\bibitem{W7}  {\sc Z. Wojtkowiak},  On l-adic Galois periods, Relations between coefficients of Galois representations on fundamental groups of a projective line minus a finite number of points, 
Actes de la conf\'erence ``Cohomologie l-adiques et corps de nombres'', 10-14 d\'ecembre 2007, CIRM Luminy,
Publ. Mathematiques de Besan\c con, Alg\`ebre et Th\'eorie des Nombres, F\'evrier 2009, pp. 157-174.
\bibitem{W8}  {\sc Z. Wojtkowiak}, On $\ell$-adic Galois L-functions, arXiv:1403.2209v1 [math. NT] 10 Mar 2014.
\end{thebibliography}
\end{document}